%% file: ElbSchSch11.tex
\documentclass[10pt]{article}
\usepackage{amsmath}
\usepackage{amssymb}
\usepackage{graphicx}
\usepackage{tikz}

\numberwithin{equation}{section}
\DeclareMathOperator{\ce}{ce}
\DeclareMathOperator{\se}{se}

\newcommand{\field}[1]{\ensuremath{\mathbb{#1}}}
\newcommand{\R}{\field{R}}
\newcommand{\N}{\field{N}}

\renewcommand{\d}{\,d}
\renewcommand{\i}{\mathrm i}
\newcommand{\e}{\mathrm e}

\newcommand{\abs}[1]{\lvert#1\rvert}
\newcommand{\tsqrt}[1]{{\textstyle\sqrt{#1}}}

\title{Reconstruction Formulas for Photoacoustic Sectional Imaging
\thanks{The work has been supported by the Austrian Science Fund (FWF)
within the national research network Photo\-acoustic Imaging in Biology and Medicine,
project S10505-N20.}}
\author{Peter Elbau$^2$, Otmar Scherzer$^{1,2}$ and Rainer Schulze$^2$
\\
\\
{\normalsize
\begin{tabular}{ccc}
\hbox to 0pt{\hss ${}^1$}Computational Science Center && \hbox to 0pt{\hss ${}^2$}Radon Institute of Computational\\
  University of Vienna &&  and Applied Mathematics\\
  Nordbergstr.~15 && Altenberger Str.~69\\
  1090 Vienna, Austria && 4040 Linz, Austria\\
\end{tabular}}}

\begin{document}

\maketitle

\begin{abstract}
The literature on reconstruction formulas for photoacoustic tomography (PAT) is vast. The various reconstruction formulas differ by
used measurement devices and geometry on which the data are sampled. In standard photoacoustic imaging (PAI), the object under
investigation is illuminated uniformly. Recently, sectional photoacoustic imaging techniques, using focusing techniques for
initializing and measuring the pressure along a plane, appeared in the literature. This paper surveys existing and provides novel
\emph{exact} reconstruction formulas for sectional photoacoustic imaging.
\end{abstract}

%%%%%%%%%%%%%%%%%%%%%%%%%%%%%%
% Introduction
%%%%%%%%%%%%%%%%%%%%%%%%%%%%%%
\section{Introduction}
The literature on reconstruction formulas and back-projection algorithms for photoacoustic imaging is vast. Wang et al.\ developed
reconstruction formulas for cylindrical, spherical, and planar measurement geometries in a series of papers~\cite{XuFenWan02,XuWan02,XuXuWan02},
and recently many more algorithms based on reconstruction formulas have been developed
(see the survey~\cite{KucKun08}).

Also different measurement devices for the ultrasound pressure have been suggested.
Most common are small detectors based on materials, which exhibit a strong piezoelectric effect and can be immersed safely in
water (i.e.\ polymers such as PVDF). In analytical reconstruction formulas, they are considered point detectors.
Other experimental setups have been realized with line and area detectors, which collect averaged pressure
(see \cite{PalNusHalBur09} for a survey).

Here, we consider the problem of photoacoustic \emph{sectional imaging}. Opposed to standard photoacoustic imaging, where the detectors
record sets of two-dimensional projection images from which the three-dimensional imaging data can be reconstructed, \emph{single slice imaging} reconstructs
a set of two-dimensional slices, each by a single scan procedure. The advantages of the latter approach are a considerable increase in measurement
speed and the possibility to do selective plane imaging. In general, this can only be obtained by the cost of decreased out-of-plane
resolution (i.e.\ the direction orthogonal to the focusing plane). Experimentally, one can obtain photoacoustic sectional imaging by
illuminating a single plane of the object and by using a focused detector. Technical details are provided in Section~\ref{sec:exp}.
In our experiments, the measurement data are recorded on a cylindrical domain $\partial \Omega \times \R$, where $\partial \Omega$ denotes
the boundary of a smooth domain $\Omega$ in $\R^2$.

The difference in this model to previously studied models is that the wave propagation is considered fully three-dimensional, the initialization
\emph{and} measurements are fully two-dimensional due to the selective plane illumination and detection. Therefore, such setups require novel
reconstruction formulas. In particular, as a further novelty, we present
reconstruction formulas in ellipsoidal domains.

This paper surveys existing and provides novel \emph{exact} formulas for the reconstruction of the initial pressure distribution for various kinds
of measurement setups. After the introduction of the \emph{universal} back-projection algorithm introduced in~\cite{XuWan05} this goal seems
superfluous, although not discussed in detail for sliced imaging. However, it has been shown recently by Natterer~\cite{Nat11}
that \emph{universal} back-projection is only exact for special sampling geometries. Here, for sliced imaging and certain sampling setups,
we can indeed find mathematically exact universal reconstruction algorithms for arbitrary strictly convex sampling domains~$\Omega$.

The paper is organized as follows:
In Section~\ref{sec:exp}, we describe the experimental setup of photoacoustic sectional imaging, and we model in Section~\ref{sec:mathematical_formulation} various measurements where it is possible to derive exact reconstructions formulas for sectional imaging.
The reconstruction formulas are then provided in Section~\ref{sec:rec}. In the appendix, we survey some background material on the Abel transform, the spherical mean operator, and the Mathieu equation.

%%%%%%%%%%%%%%%%%%%%%%%%%%%%%%
% Experimental Background
%%%%%%%%%%%%%%%%%%%%%%%%%%%%%%
\section{Experimental Background}
\label{sec:exp}
Below we give an overview on \emph{photoacoustic sectional imaging}, describe the experimental realization, and provide mathematical
formulations. Opposed, we call conventional photoacoustic imaging with uniform illumination of the object \emph{non-focused}.

In general, PAI is based on the so called thermo- or photoacoustic effect. Laser light impinging onto a surface of an object leads to its heating and under the conditions of thermal and stress confinement an acoustic wave emerges from the object. Mathematically, this means that the initial pressure distribution $f$ is related to the absorption coefficient of the object by~\cite{CoxArrBea06}
\begin{equation}
\label{eq:rel_pr_abs}
f(x) = \frac{\beta c^2}{C_p}\mu_a(x)\Phi(x)
\end{equation}
where $x\in R^3$. The constant prefactor consists of $\beta$ the thermal expansivity, $C_p$ the specific heat capacity and $c$ the speed of sound. Then, $\mu_a$ is the absorption coefficient of the object and $\Phi$ the local light fluency. Neglecting effects of light propagation (i.e. setting $\Phi$ a constant), $\mu_a$ remains the only variable depending on $x$ and containing the desired tomographic information about the object.
Reconstruction algorithms that also take light propagation into account are investigated in Ref.~\cite{CoxArrBea06} but are not focus of this work.

It is common to classify photoacoustic measurement setups into \emph{point} (see e.g.~\cite{KruReiKru99}) and \emph{integrating detector}
setups~\cite{BurHofPalHalSch05}. The focus of this paper are reconstruction methods for sectional imaging and various kinds of detectors.

\begin{figure}[htb]
\centerline{\includegraphics[width=10cm]{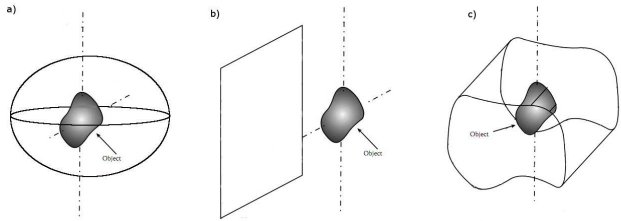}}
\caption{Sketch of non-focused point detector arrays. a) closed surface, b) quasi-infinite plane, c) cylindrical surface. The point detectors are spread over the indicated surfaces (points are not drawn in the figure). The detector arrays b) and c) suffer from the limited view problem.}
\label{fig:pointgeom}
\end{figure}

In conventional photoacoustic setups with point detectors, measurements are collected (see Fig.~\ref{fig:pointgeom})
all over a surface enclosing the object (i.e.\ a sphere or an ellipsoid), or all over a sufficiently large cylindrical surface
(where the cross-section can be every line segment), or over a sufficiently large plane.
The latter two detector array geometries allow for approximate reconstructions only, since sufficiently large means that in
theory the array is considered infinitely large. This approximation causes the limited view problem in practical applications.
Paltauf et al.\ have given correction factors for some of the affected geometries~\cite{PalNusBur09}.

The linear and planar detectors have to be moved tangentially to a surface surrounding the object (see~\cite{BurHofPalHalSch05}).
Practically, this only allows the measurement devices to be aligned on a cylindrical surface (or on a plane).
Experimental realizations of line detectors are documented for instance in~\cite{PalNusHalBur09}.
\begin{figure}[htb]
\centerline{\includegraphics[width=10cm]{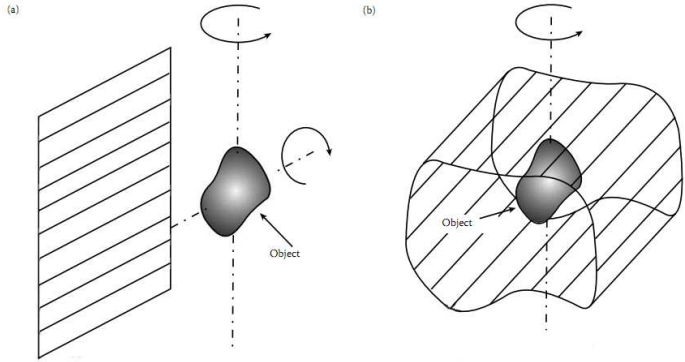}}
\caption{Sketch of non-focused line detector arrays. a) quasi-infinite plane, b) cylindrical surface with arbitrary cross-section. The detector array a) suffers from the limited view problem.}
\label{fig:linegeom}
\end{figure}

In the following, we explain the principles of focusing detectors. The ultrasonic wave is refracted by a suitable acoustic lens such that
out-of-plane signals are generally weak and can be neglected. Thus only signals emerging from the desired imaging plane are
collected at the detector.
Contemporary focusing ultrasonic detectors have a spherical or cylindrical shape, thus the detector surface plays the role of the acoustic
lens. We consider the case of a cylindrically focusing detector, which focuses into a plane.
The sectional imaging can be improved further by illuminating the desired plane only, i.e.\ by cylindrical lenses in front of the object. Note, however, that this requires a low scattering coefficient of the sample, because otherwise illumination will not be restricted to a single plane. However, both approaches in combination provide a good experimental approximation to the mathematical models outlined below. The out-of-plane resolutions achieved are worse than the in-plane resolutions nevertheless. For more details on focusing point detectors see~\cite{MaTarNtzRaz09,RazDisVinMaPerKosNtz09} and for focusing line detectors see~\cite{GraPasNusPal10,GraPasNusPal11}.

%%%%%%%%%%%%%%%%%%%%%%%%%%%%%%
% Mathematical Formulation
%%%%%%%%%%%%%%%%%%%%%%%%%%%%%%
\section{Mathematical Formulation}\label{sec:mathematical_formulation}
We model the sectional photoacoustic imaging by assuming that the \emph{initial pressure distribution} $f:\R^3\to\R$ is perfectly focused in the illumination plane $\{x\in\R^3\mid x_3=0\}$:
\begin{equation}
\label{eq:hatf}
f(\xi,z) = \hat f(\xi) \delta(z),\quad\xi\in\R^2,\;z\in\R,
\end{equation}
for some smooth function $\hat f:\R^2\to\R$. The resulting pressure wave $p:\R^3\times[0,\infty)\to\R$, we assume to be the solution of the linear three-dimensional wave equation
\begin{equation}
\label{eq:wave}
\boxed{
\begin{aligned}
\partial_{tt}p(\xi,z;t) &= \Delta_{\xi,z} p(\xi,z;t),\\
\partial_tp(\xi,z;0) &= 0, \phantom{\hat f}\\
p(\xi,z;0) &= f(\xi,z) = \hat f(\xi) \delta(z)
\end{aligned}}
\end{equation}
for all $\xi\in\R^2$, $z\in\R$, and $t>0$. Here,
\begin{equation*}
\Delta_{\xi,z} = \partial_{\xi_1\xi_1}+\partial_{\xi_2\xi_2}+\partial_{zz}
\end{equation*}
denotes the three-dimensional Laplacian in Euclidean coordinates.

Our aim is to recover the function $\hat f$, describing the initial pressure distribution, from certain messurements of the pressure wave $p$. The position of the detectors performing these measurements shall be given by the boundary $\partial\Omega$ of a convex domain $\Omega\subset\R^2$ in the illumination plane, where we additionally assume that $\hat f$ has compact support in $\Omega$.

We will consider the following four different measurement setups and derive reconstruction formulas. Some of the reconstruction formulas and setups have already been documented in the literature and are surveyed here in the general context. Most of the reconstruction formulas, however, are new.

\begin{description}
\item[Vertical Line Detectors:]
The measurement data are
\begin{equation}\label{eq:m1}
\boxed{
\begin{aligned}
m_1(\xi;t) &:=  \int_{-\infty}^\infty p(\xi,z;t)\d z \\
& \text{ for all } \xi \in \partial\Omega,\;t>0.
\end{aligned}}
\end{equation}
That is, in practical realizations, we use \emph{line detectors} which measure the overall pressure along a line orthogonal to the illumination plane.

\item[Point Detectors:]
The measurement data are
\begin{equation}\label{eq:m2}
\boxed{
\begin{aligned}
m_2(\xi;t) & :=  p(\xi,0;t) \\
& \text{ for all } \xi \in\partial\Omega,\;t>0.
\end{aligned}}
\end{equation}
That is, we use \emph{point detectors} which measure the pressure on the boundary of $\partial\Omega$ over time. This measurement geometry is used in~\cite{MaTarNtzRaz09,RazDisVinMaPerKosNtz09}.
\end{description}

For the other two measurement methods, we additionally impose that the domain $\Omega\subset\R^2$ is strictly convex and bounded.
\begin{description}
\item[Vertical Plane Detectors:]
The measurement data are
\begin{equation}\label{eq:m3}
\boxed{
\begin{aligned}
m_3(\theta;t) &:=  \int_{P(\theta)} p(x;t)\d s(x) \\
& \text{ for all } \theta \in S^1,\;t>0,
\end{aligned}}
\end{equation}
where $P(\theta)\subset\R^3$ denotes the tangential plane of the cylinder $\partial\Omega\times\R$ orthogonal to the vector~$(\theta,0)$, see~\eqref{eq:tangent_planes}. For a practical realization we use \emph{planar detectors} which are moved tangentially around the object and measure the averaged pressure on the plane.

\item[Horizontal Line Detectors:]
The measurement data are
\begin{equation}\label{eq:m4}
\boxed{
\begin{aligned}
m_4(\theta;t) &:=  \int_{T(\theta)} p(\xi,0;t)\d s(\xi)\\
& \text{ for all }  \theta \in S^1,\;t>0,
\end{aligned}}
\end{equation}
where $T(\theta)\subset\R^2$ denotes the tangential line of $\partial\Omega$ orthogonal to the vector~$\theta$, see~\eqref{eq:tangent_lines}. This is a realization using \emph{line detectors} which measure the overall pressure on a line tangential to $\partial\Omega$ in the illumination plane, see \cite{GraPasNusPal10,GraPasNusPal11}.
\end{description}

In those cases where the domain $\Omega$ is strictly convex and bounded, we parametrize the boundary $\partial\Omega$ with the map $\zeta:S^1\to\partial\Omega$ which associates to every unit vector $\theta\in S^1$ the point $\zeta(\theta)\in\partial\Omega$ where the outward unit normal vector of $\partial\Omega$ coincides with $\theta$, see Figure~\ref{fig:zeta}.
\input{point_of_tangency.tex}

Since the tangent line $T(\theta)$ of $\partial\Omega$ at $\zeta(\theta)$ is thus by definition orthogonal to~$\theta$, we can define the family $T(r,\theta)$, $r\in\R$, of lines parallel to the tangent~$T(\theta)$ by
\begin{equation}\label{eq:tangent_lines}
T(r,\theta)= \zeta(\theta)+r\theta+\R\theta^\perp\subset \R^2,\quad T(\theta)=T(0,\theta),
\end{equation}
for every $\theta\in S^1$ and $r\in\R$.

Moreover, we introduce the family $P(r,\theta)$, $r\in\R$, of planes parallel to the tangent plane $P(\theta)$ of the cylinder $\partial\Omega\times\R$ at $(\zeta(\theta),0)$ by
\begin{equation}
\label{eq:tangent_planes}
P(r,\theta)=(T(\theta),0)+(0,\R) \subset \R^3,\quad P(\theta)=P(0,\theta),
\end{equation}
for every $\theta\in S^1$ and $r\in\R$.

%%%%%%%%%%%%%%%%%%%%%%%%%%%%%%
% Reconstruction Methods
%%%%%%%%%%%%%%%%%%%%%%%%%%%%%%
\section{Reconstruction Methods}
\label{sec:rec}
In the following, we derive universal reconstruction formulas for photoacoustic sectional imaging. Conceptually, the paper is closely
related to~\cite{Kun07b}, where universal formulas for conventional photoacoustic tomography were derived in arbitrary geometry. However,
some of the results there are implicit, requiring explicit knowledge of eigenfunctions, which are provided here explicitly.
Even more, the focus of this paper is on sectional imaging, which results in different formulas.

%%%%%%%%%%%%%%%%%%%%%%%%%%%%%%
% m_1
%%%%%%%%%%%%%%%%%%%%%%%%%%%%%%
\subsection{Measurements with Vertical Line Detectors}
We introduce the function
\begin{equation}\label{eq:line_o}
\tilde p(\xi;t)=\int_{-\infty}^\infty p(\xi,z;t)\d z,\quad \xi\in\R^2,\;t\ge0.
\end{equation}

Then the inital value problem~\eqref{eq:wave} for the function $p$ implies that the function $\tilde p$ satisfies the two-dimensional wave equation
\[ \partial_{tt}\tilde p(\xi;t) = \Delta_\xi\tilde p(\xi;t)\quad\text{for all}\quad\xi\in\R^2,\;t>0 \]
with the initial conditions
\begin{alignat*}{2}
\partial_t\tilde p(\xi;0)&=0&&\text{for all}\quad\xi\in\R^2, \\
\tilde p(\xi;0)&=\hat f(\xi)&\quad&\text{for all}\quad\xi\in\R^2.
\end{alignat*}
The initially three-dimensional reconstruction problem therefore reduces to the two-dimensional problem of calculating $\hat f(\xi)=\tilde p(\xi;0)$, $\xi\in\R^2$, from the measurement data
\[ m_1(\xi;t) = \tilde p(\xi;t),\quad \xi\in\partial\Omega,\;t>0. \]

\subsubsection{Reconstruction Formulas Based on Series Expansions}
For special domains $\Omega$, explicit reconstruction formulas are known: see the review \cite{HalSchBurNusPal07} for $\Omega$ a
circle and the half-space. The derivation for the ellipse is published in \cite{ElbSchSch12}.

\begin{itemize}
\item
If $\Omega$ is the half-space $\{\xi\in\R^2\mid\xi_2>0\}$, we get~\cite{KosBea03}
\begin{equation}\label{eq:NorLin1}
\boxed{
\hat f(\xi) = \frac2\pi\int_{-\infty}^\infty\int_k^\infty\tilde m_1(k,\omega)\e^{\i k\xi_1}\cos(\xi_2\tsqrt{\omega^2-k^2})\d\omega\d k
}
\end{equation}
for every $\xi\in\R^2$, where
\[ \tilde m_1(k,\omega) = \frac1\pi\int_{-\infty}^\infty\int_0^\infty m_1(\xi_1,0;t)\e^{-\i k\xi_1}\cos(\omega t)\d t\d\xi_1 \]
is the Fourier--cosine transform of the measurement data $m_1$.
\item
If $\Omega=B^2_R(0)$ is the two-dimensional ball with radius $R$ and center $0$, we choose $\psi:[0,\infty)\times[0,2\pi)\to\R^2$, $\psi(r,\varphi)=(r\cos\varphi,r\sin\varphi)$, and get~\cite{HalSchBurNusPal07}
\begin{equation}\label{eq:Nor1}
\boxed{
\hat f(\psi(r,\varphi))=\frac1\pi\int_0^\infty\sum_{k=-\infty}^\infty\frac{J_{|k|}(R\omega r)}{J_{|k|}(R\omega)}\,\tilde m_1(k,\omega)\e^{\i k\varphi}\d\omega
}
\end{equation}
for every $r\in[0,\infty)$ and $\varphi\in[0,2\pi)$, where
\[ \tilde m_1(k,\omega)=\frac1\pi\int_0^\infty\int_0^{2\pi}m_1(\psi(R,\varphi);t)\e^{-\i k\varphi}\cos(\omega t)\d\varphi\d t \]
is the Fourier--cosine transform of the measurement data $m_1$.
Here, $J_k$, $k\in\N_0$, denotes the $k$th Bessel function.
\end{itemize}

Equations~\eqref{eq:Nor1} and~\eqref{eq:NorLin1} can be derived from formulas for the inversion of the spherical mean operator (see Section~\ref{sec:spherical_mean}), and this is why these formulas are typically assigned to Norton~\cite{Nor80} and Norton \& Linzer~\cite{NorLin81}, although they considered reflectivity ultrasound imaging and here the topic is photoacoustics.
\begin{itemize}
\item
If $\Omega$ is the ellipse $\{\xi\in\R^2\mid\tfrac{\xi_1^2}{a^2}+\tfrac{\xi_2^2}{b^2}<1\}$ with $a>b$, we set
\[ \psi:[0,\infty)\times[0,2\pi)\to\R^2,\quad\psi(r,\varphi)=\varepsilon\begin{pmatrix}\cosh(r)\cos(\varphi)\\\sinh(r)\sin(\varphi))\end{pmatrix},\quad \]
with the linear eccentricity $\varepsilon=\sqrt{a^2-b^2}$ and find~\cite{ElbSchSch12}
\[ \boxed{ \hat f(\psi(r,\varphi))=\frac{\sqrt2}\pi\int_0^\infty\sum_{k=0}^\infty\frac{R_k(r;\omega)}{R_k(r_0;\omega)}\tilde m_1(k,\omega)\Phi_k(\varphi;\omega)\d\omega, } \]
where $r_0=\mathrm{artanh}(\frac ba)$ is chosen such that $\psi(r_0,\varphi)\in\partial\Omega$ and
\[ \tilde m_1(k,\omega) = \frac{\sqrt2}\pi\int_0^\infty\int_0^{2\pi}m_1(\psi(r_0,\varphi);t)\Phi_k(\varphi;\omega)\cos(\omega t)\d\varphi\d t. \]
Herein, the functions $\Phi_k$ are for $k\in\N_0$ defined by
\[ \Phi_{2k}(\varphi;\omega)=\ce_k(\varphi;\tfrac{\varepsilon^2\omega^2}4)\quad\text{and}\quad\Phi_{2k+1}(\varphi;\omega)=\se_{k+1}(\varphi;\tfrac{\varepsilon^2\omega^2}4), \]
where $\ce_k$ and $\se_k$ denote the Mathieu cosine and Mathieu sine functions, respectively, see Section~\ref{sec:Mathieu_equation}, and the functions $R_k$ are the corresponding solutions of the radial Mathieu equation and are for $k\in\N_0$ given by
\[ R_{2k}(r;\omega)=\ce_k(\i r;\tfrac{\varepsilon^2\omega^2}4)\quad\text{and}\quad R_{2k+1}(r;\omega)=-\i\se_{k+1}(\i r;\tfrac{\varepsilon^2\omega^2}4). \]
\end{itemize}

\subsubsection{Reduction to the Spherical Mean Operator}
\label{se:reductionToSphericalMean}
Taking into account the relation~\eqref{eq:spmeans_and_2d_wave_equation} between the solution of the two-dimen\-sional wave equation and the spherical mean operator, the problem of determining $\hat f$ from the measurements $m_1(\xi;t)=\tilde p(\xi;t)$ for $\xi\in\partial\Omega$, $t>0$ can be equivalently described as the problem of reconstructing $\hat f$ from the spherical mean operator $\mathcal M_2[\hat f](\xi;r)$ of $\hat f$ for $\xi\in\partial\Omega$ and $r\in(0,\infty)$.

\begin{itemize}
\item
For $\Omega = B^2_R(0) \subset \R^2$, analytical reconstruction formulas have been derived by Finch, Haltmeier, Rakesh \cite{FinHalRak07} and read as follows
\begin{equation}
\label{m-lap}
\boxed{
\hat{f}(\xi) = \frac{1}{2 \pi} \Delta_{\xi} \left(\int_{S^1} \int_0^{2R} r (\mathcal M_2[\hat{f}])(R\theta,r) \log \abs{r^2 -\abs{ \xi - R\theta}^2 } \d r \d s(\theta)\right)}
\end{equation}
and
\begin{equation}\label{m-inv}
\boxed{
\hat f(\xi) = \frac{1}{2\pi} \int_{S^1} \int_0^{2R}\left( \partial_r r \partial_r \mathcal M_2[\hat{f}] \right)(R\theta,r) \log \abs{r^2 - \abs{\xi - R\theta}^2}\d r\d s(\theta).}
\end{equation}
\item
For a general domain~$\Omega$, Kunyansky reduced in~\cite{Kun07b} the reconstruction problem to the determination of the eigenvalues $\lambda_k$ and normalized eigenfunctions $u_k$, $\|u_k\|_2=1$, of the Dirichlet Laplacian $-\Delta$ on $\Omega$ with zero boundary conditions:
\begin{align}
\Delta u_k(\xi)+\lambda_ku_k(\xi)&=0,\quad\xi\in\Omega, \label{eq:Helmholtz_equation}\\
u_k(\xi)&=0,\quad\xi\in\partial\Omega.
\end{align}
Indeed, if $(\xi,\eta)\mapsto G_{\lambda_k}(|\xi-\eta|)$ is a free-space rotationally invariant Green's function of the Helmholtz equation~\eqref{eq:Helmholtz_equation} and $n(\xi)$ denotes the outer unit normal vector of $\partial\Omega$ at $\xi\in\partial\Omega$, then
\begin{equation}\label{eq:kunyansky}
\boxed{\hat f(\xi) = 2\pi\sum_{k=0}^\infty\tilde M_ku_k(\xi),}
\end{equation}
where
\[ \tilde M_k = \int_{\partial\Omega}\int_0^\infty r\mathcal M_2[\hat f](\eta,r)G_{\lambda_k}(r)\left<\nabla u_k(\eta),n(\eta)\right>\d r\d s(\eta). \]
\end{itemize}

We thus get the initial pressure distribution $\hat f$ by first calculating from the measurements $m_1(\xi;t)=\tilde p(\xi;t)$, $\xi\in\partial\Omega$, $t>0$, with formula~\eqref{eq:spmeans_and_2d_wave_equation} the spherical mean operator $\mathcal M_2[\hat f](\xi;r)$ of $\hat f$ for $\xi\in\partial\Omega$ and $r\in(0,\infty)$, and then using one of the formulas \eqref{m-lap}, \eqref{m-inv}, or~\eqref{eq:kunyansky} to get~$\hat f$.

%%%%%%%%%%%%%%%%%%%%%%%%%%%%%%
% m_2
%%%%%%%%%%%%%%%%%%%%%%%%%%%%%%
\subsection{Measurements with Point Detectors}
From equation~\eqref{ndsol3}, we know that the solution of the initial value problem~\eqref{eq:wave} can be for every $x\in\R^3$ and $t>0$ written in the form
\[ p(x;t) = \partial_t\left(\frac1{4\pi t}\int_{\partial B^3_t(0)}f(x+y)\d s(y)\right). \]
Parameterizing the sphere $\partial B^3_t(0)$ in cylindrical coordinates, i.e.\ in the form
\[ \partial B^3_t(0)=\left\{(\tsqrt{t^2-h^2}\,\theta,h)\;:\;h\in[-t,t],\;\theta\in S^1\right\}, \]
we find for every $x=(\xi,z)$, $\xi\in\R^2$, $z\in\R$, and $t>0$ that
\[ p(\xi,z;t) = \partial_t\left(\frac1{4\pi t}\int_{-t}^t\int_{S^1}\hat f(\xi+\tsqrt{t^2-h^2}\,\theta)\delta(z+h)t\d s(\theta)\d h\right), \]
where we have used the special property of the source term \eqref{eq:hatf}.
Integrating out the $\delta$-distribution, we get for $z\in[-t,t]$
\begin{equation}
\label{eq:m2_central}
p(\xi,z;t)=\partial_t\left(\frac1{4\pi}\int_{S^1}\hat f(\xi+\tsqrt{t^2-z^2}\,\theta)\d s(\theta)\right).
\end{equation}
By the definition~\eqref{spmean} of the spherical mean operator $\mathcal M_2$, this means
\begin{equation}\label{eq:m2_p_as_derivative}
p(\xi,z;t)=\frac12\partial_t\big(\mathcal M_2[\hat f](\xi;\tsqrt{t^2-z^2})\big)\quad\text{for}\quad z\in[-t,t].
\end{equation}
From the assumption that the support of $\hat{f}$ lies completely in $\Omega$, we know that $\mathcal M_2[\hat f](\xi;0)=\hat f(\xi)=0$ for $\xi\notin\Omega$. Thus, we can integrate the relation~\eqref{eq:m2_p_as_derivative} for $\xi\notin\Omega$ and find for every $z\in[-t,t]$ that
\[ \mathcal M_2[\hat f](\xi;\tsqrt{t^2-z^2}) = 2\int_z^tp(\xi,z;\tilde t)\d\tilde t. \]
Setting $z=0$, we get for every $\xi\in\partial\Omega$ and every $t>0$ the relation
\[ \boxed{\mathcal M_2[\hat f](\xi;t)=2\int_0^tm_2(\xi;\tilde t)\d\tilde t.} \]
Having calculated the spherical mean of $\hat f$, we can now proceed as in Section~\ref{se:reductionToSphericalMean}.

%%%%%%%%%%%%%%%%%%%%%%%%%%%%%%
% m_3
%%%%%%%%%%%%%%%%%%%%%%%%%%%%%%
\subsection{Measurements with Vertical Plane Detectors}\label{sec_m3}
For every $\theta\in S^1$, we define for $r\in\R$ and $t\ge0$ the function
\[ \tilde p_\theta(r;t) = \int_{P(r,\theta)}p(x;t)\d s(x)\,,\]
where $P(r,\theta)$ denotes the plane as defined in~\eqref{eq:tangent_planes}.

Then, since the vectors $(\theta,0)$, $(\theta^\perp,0)$, and $(0,0,1)$ form an orthonormal basis of~$\R^3$ and the Laplacian is rotationally invariant, we find from equation~\eqref{eq:wave} that
\[ \partial_{tt}\tilde p_\theta(r;t)=\int_{-\infty}^\infty\int_{-\infty}^\infty\Delta_x p(\zeta(\theta)+r\theta+u\theta^\perp,z;t)\d u\d z = \partial_{rr}\tilde p_\theta(r;t) \]
for every $r\in\R$ and $t>0$.
Thus, $\tilde p_\theta$ solves the one-dimensional wave equation with the initial conditions
\begin{alignat*}{2}
\partial_t\tilde p_\theta(r;0)&=0&&\text{for all $r\in\R$ and} \\
\tilde p_\theta(0;t)&=m_3(\theta;t)&\quad&\text{for all $t>0$}
\end{alignat*}
resulting from~\eqref{eq:wave} and~\eqref{eq:m3}, respectively.
Moreover, since $\hat f$ has its support inside~$\Omega$, we know that $\tilde p_\theta(r;0)=0$ for $r\ge0$.

With d'Alembert's formula for the solution of the one-dimensional wave equation, we find that the unique solution for this initial value problem is given by
\[ \tilde p_\theta(r;t) = m_3(\theta;-t-r)+m_3(\theta;t-r),\quad r\in\R,\;t>0, \]
where we set $m_3(\theta;t)=0$ for $t\le0$.

Finally, we have to recover from the values of $\tilde p_\theta$, $\theta\in S^1$, the initial pressure distribution $\hat{f}$ from equation~\eqref{eq:wave}. We have the relation
\[ \tilde p_\theta(r;0)=\int_{-\infty}^\infty \hat{f}(\zeta(\theta)+r\theta+u\theta^\perp)\d u = \mathcal R[\hat f](r+\left<\zeta(\theta),\theta\right>,\theta), \]
where $\mathcal R$ denotes the Radon transform as defined in~\eqref{eq:radon_transform}. We can therefore recover $\hat f$ with an inverse Radon transform:
\begin{equation}
\label{eq:ir}
\boxed{\hat f = 2\mathcal R^{-1}[\tilde m_3],\quad \tilde m_3(r,\theta)=\begin{cases}m_3(\theta;\left<\zeta(\theta),\theta\right>-r)&\text{if}\;r<\left<\zeta(\theta),\theta\right>,\\\;0&\text{if}\;r\ge\left<\zeta(\theta),\theta\right>.\end{cases}}
\end{equation}
Equation~\eqref{eq:ir} reveals an interesting property of integrating area detectors: For an arbitrary strictly convex measurement geometry
$\Omega$, exact reconstruction formulas exist. This is a property which is not known for conventional and other
photoacoustic sectional imaging technologies.

%%%%%%%%%%%%%%%%%%%%%%%%%%%%%%
% m_4
%%%%%%%%%%%%%%%%%%%%%%%%%%%%%%
\subsection{Measurements with Horizontal Line Detectors}
For every $\theta\in S^1$, we define the function
\[ \tilde p_\theta(r,z;t)=\int_{T(r,\theta)}p(\xi,z;t)\d s(\xi), \]
where $T(r,\theta)$ is defined as in~\eqref{eq:tangent_lines}. Then, using that the vectors $(\theta,0)$, $(\theta^\perp,0)$, and $(0,0,1)$ are an orthonormal basis of $\R^3$ and that the Laplacian is rotationally invariant, the initial value problem~\eqref{eq:wave} implies that $\tilde p_\theta$ solves for all $r,z\in\R$ and $t>0$ the two-dimensional wave equation
\begin{align*}
\partial_{tt}\tilde p_\theta(r,z;t) &= \int_{-\infty}^\infty\Delta_x p(\zeta(\theta)+r\theta+u\theta^\perp,z;t)\d u \\
&=\partial_{rr}\tilde p_\theta(r,z;t)+\partial_{zz}\tilde p_\theta(r,z;t)
\end{align*}
with the initial conditions
\begin{align*}
\partial_t\tilde p_\theta(r,z;0)&=0, \\
\tilde p_\theta(r,z;0)&=F_\theta(r)\delta(z),\quad F_\theta(r)=\int_{T(r,\theta)}\hat f(\xi)\d s(\xi),
\end{align*}
for every $r,z\in\R$.

From formula~\eqref{eq:2d_wave_equation_solution}, we see that the solution of this initial value problem can be written as
\begin{align*}
\tilde p_\theta(r,z;t) &= \frac1{2\pi}\partial_t \left(\int_{B^2_t(0)}\frac{F_\theta(r+\rho)\delta(z+\zeta)}{\sqrt{t^2-\rho^2-\zeta^2}}\d s(\rho,\zeta)\right) \\
&= \frac1{2\pi}\partial_t\left(\int_{-t}^t\delta(z+\zeta)\int_{-\sqrt{t^2-\zeta^2}}^{\sqrt{t^2-\zeta^2}}\frac{F_\theta(r+\rho)}{\sqrt{t^2-\rho^2-\zeta^2}}\d\rho\d\zeta\right)
\end{align*}
for all $r,z\in\R$ and $t>0$.
Integrating out the $\delta$-function, we find for every $z\in[-t,t]$ that
\[ \tilde p_\theta(r,z;t) = \frac1{2\pi}\partial_t \left(\int_{-\sqrt{t^2-z^2}}^{\sqrt{t^2-z^2}}\frac{F_\theta(r+\rho)}{\sqrt{t^2-z^2-\rho^2}}\d\rho\right). \]

Since $\tilde p_\theta$ is related to the measurement $m_4$, given by~\eqref{eq:m4}, via $m_4(\theta;t) = \tilde p_\theta(0,0;t)$, and since $F_\theta(r)=0$ for $r>0$ by the assumption that $\hat f$ has support inside $\Omega$, we find with the formula~\eqref{eq:abel} for the Abel transform in reciprocal coordinates that
\[ m_4(\theta;t) = \frac1{2\pi}\partial_t \left(\int_0^t\frac{F_\theta(-\rho)}{\sqrt{t^2-\rho^2}}\d\rho\right) = \frac1{4\pi}\partial_t\left(\tfrac1t\mathcal A[\psi_\theta](\tfrac1t)\right) \]
where $\psi_\theta(\tfrac1{\rho})=\rho^2F_\theta(-\rho)$. Switching to the reciprocal coordinate $s=\frac1t$ and using the identity~\eqref{eq:inverse_abel}, we see that this is of the form
\[ \frac2{s^2}m(\theta;\tfrac1s) = -\frac1{2\pi}\partial_s\big(s\mathcal A[\psi_\theta](s)\big) = \mathcal A^{-1}[\tilde\psi_\theta](s) \]
with $\tilde\psi_\theta(\frac1\rho)=\frac1{\rho^2}\psi_\theta(\frac1\rho)=F_\theta(-\rho)$.
Thus, we can directly solve the equation for $F_\theta$ and find
\begin{equation}\label{eq:m4_formula_for_F_theta}
F_\theta(-\rho) = 2\mathcal A[\tilde m_\theta](\tfrac1\rho),\quad\tilde m_\theta(\tfrac1t) = t^2m_4(\theta;t).
\end{equation}

Since we have by definition
\[ F_\theta(r) = \mathcal R[\hat f](r+\left<\zeta(\theta),\theta\right>,\theta), \]
we finally get (remembering that $F_\theta(r)=0$ for $r\ge0$)
\[ \boxed{\hat f=2\mathcal R^{-1}[\tilde F],\quad \tilde F(r,\theta)=\begin{cases}\mathcal A[\tilde m_\theta]\Big(\tfrac1{\left<\zeta(\theta),\theta\right>-r}\Big)&\text{if}\;r<\left<\zeta(\theta),\theta\right>,\\\;0&\text{if}\;r\ge\left<\zeta(\theta),\theta\right>.\end{cases}} \]
So, the reconstruction of $\hat f$ can be accomplished by an Abel transform of the rescaled measurements $\tilde m_\theta$, defined in~\eqref{eq:m4_formula_for_F_theta}, followed by an inverse Radon transform. Again, this reconstruction formula is valid for an arbitrary strictly convex measurement geometry~$\Omega$.

\section*{Conclusion}
In this paper we have surveyed exact reconstruction formulas for photoacoustic sectional imaging. All formulas are mathematically, analytically
exact. Comparing point and integrating line detectors, it is quite surprising that integrating area detectors allow analytical reconstructions for
all strictly convex domains. This property has not been observed for point detectors.

\appendix
\section{Appendix}
\subsection{Abel Transform}
The \emph{Abel transform} $\mathcal A[\psi]$ of a smooth function $\psi: \R_+ \to \R$, which decays to zero at $\infty$, is defined by
\begin{equation*}
\mathcal A[\psi](y) = \int_{-\infty}^\infty\psi(\tsqrt{x^2+y^2})\d x = 2\int_y^\infty \frac{r\psi(r)}{\sqrt{r^2-y^2}}\d r, \quad y\ge0.
\end{equation*}
We rewrite the Abel transform in reciprocal coordinates, so that it better fits in our context. Substituting $y=\frac1t$ and $r=\frac1s$, we find that
\begin{equation}\label{eq:abel}
\mathcal A[\psi](\tfrac1t) = 2\int_{\frac1t}^\infty\frac{tr\psi(r)}{\sqrt{r^2t^2-1}}\d r = 2\int_0^t\frac{t\psi(\tfrac1s)}{s^2\sqrt{t^2-s^2}}\d s.
\end{equation}

To invert the Abel transform, we remark that we have for all $v\ge0$
\begin{align*}
\int_{-\infty}^\infty\frac{(\mathcal A[\psi])'(\sqrt{u^2+v^2})}{\sqrt{u^2+v^2}}\d u &= \int_{-\infty}^\infty\int_{-\infty}^\infty\frac{\psi'(\sqrt{x^2+u^2+v^2})}{\sqrt{x^2+u^2+v^2}}\d u\d x \\
&= 2\pi\int_v^\infty\psi'(\rho)\d\rho = -2\pi\psi(v),
\end{align*}
where we substituted $x=\sqrt{\rho^2-v^2}\cos(\varphi)$ and $u=\sqrt{\rho^2-v^2}\sin(\varphi)$.
Therefore, the \emph{inverse Abel transform} $\mathcal A^{-1}[\psi]$ of a function $\psi:\R_+\to\R$ can be written as
\begin{equation}
\label{eq:diff_property}
 \mathcal A^{-1}[\psi](y) = -\frac1{2\pi y}(\mathcal A[\psi])'(y).
 \end{equation}

Using the identity
\begin{align*}
(\mathcal A[r^2\psi])'(y) &= \partial_y\left(\int_{-\infty}^\infty(x^2+y^2)\psi(\tsqrt{x^2+y^2})\d x\right) \\
&=\partial_y\big(y^2\mathcal A[\psi](y)\big)+\int_{-\infty}^\infty\frac{x^2y}{\sqrt{x^2+y^2}}\psi'(\tsqrt{x^2+y^2})\d x \\
&= \partial_y\big(y^2\mathcal A[\psi](y)\big)-y\int_{-\infty}^\infty\psi(\tsqrt{x^2+y^2})\d x = y\partial_y\big(y\mathcal A[\psi](y)\big),
\end{align*}
and using (\ref{eq:diff_property}) we can also write the inverse Abel transform in the form
\begin{equation}\label{eq:inverse_abel}
\mathcal A^{-1}[r^2\psi](y) = -\frac1{2\pi}\partial_y\big(y\mathcal A[\psi](y)\big).
\end{equation}

\subsection{Circular and Spherical Means}\label{sec:spherical_mean}
Let $x \in \R^n$ and $r\ge0$. The \emph{spherical mean operator} in $\R^n$ of an integrable function $f:\R^n\to\R$ is defined by
\begin{equation} \label{spmean}
\mathcal M_n[f](x;r) = \frac1{\abs{S^{n-1}}} \int_{S^{n-1}} f(x + r \theta)\d s(\theta),
\end{equation}
where $\abs{S^{n-1}}$ denotes the area of the unit sphere $S^{n-1}$ in $\R^n$.

The spherical mean value operator is closely related to the solution of the $n$-dimensional wave equation
\begin{align*}
\partial_{tt}p(x;t)&=\Delta_xp(x;t), \\
\partial_tp(x;0)&=0, \\
p(x;0)&=f(x)
\end{align*}
for all $x\in\R^n$ and $t>0$. More precisely, the solution $p$ can be expressed in terms of the spherical mean operator of $f$ by (see e.g.\ \cite{Eva98})
\begin{equation*}%\label{ndsol}
p(x;t) = \frac{1}{(n-2)!}\partial_t^{n-1}\left(\int_0^t r (t^2 - r^2)^{(n-3)/2} \mathcal M_n[f](x;r)\d r\right)
\end{equation*}
for all $x\in\R^n$ and $t>0$.
In particular, we have for
\begin{itemize}
\item $n=2$ that the solution $p$ of the two-dimensional wave equation can be calculated from the spherical means via the Abel transform~\eqref{eq:abel}:
\begin{equation}\label{ndsol2}
p(x;t) = \partial_t\left(\int_0^t \frac{r \mathcal M_2[f](x;r)}{\sqrt{t^2 - r^2}}\d r\right)= \partial_t\left(\frac1{2t}\mathcal A[\tilde f_x](\tfrac1t)\right),
\end{equation}
where $\tilde f_x(\tfrac1r):=r^3\mathcal M_2[f](x,r)$, leading also to the formula
\begin{equation}\label{eq:2d_wave_equation_solution}
p(x;t) = \frac1{2\pi}\partial_t\left(\int_{B^2_t(0)}\frac{f(x+y)}{\sqrt{t^2-|y|^2}}\d s(y)\right)
\end{equation}
for all $x\in\R^2$, $t>0$, where $B^2_t(0)\subset\R^2$ denotes the two-dimensional ball with radius $t$ and center $0$;

\item and for $n=3$, we get that the solution $p$ of the three-dimensional wave equation and $\mathcal M_3[f]$ are related by
\begin{equation}\label{ndsol3}
p(x;t) = \partial_t\big(t\mathcal M_3[f](x;t)\big) = \partial_t\left(\frac1{4\pi t}\int_{\partial B^3_t(0)}f(x+y)\d s(y)\right)
\end{equation}
for all $x\in\R^3$, $t>0$, where $\partial B^3_t(0)$ denotes the boundary of the three-dimensional ball $B^3_t(0)\subset\R^3$ with radius $t$ and center $0$.
\end{itemize}

We remark that we can solve the equations~\eqref{ndsol2} and~\eqref{ndsol3} for the spherical mean operator of $f$. We get for $n=2$ with $s=\frac1t$ that
\[ \frac1{\pi s^2}p(x;\tfrac1s) = -\frac1{2\pi}\partial_s\big(s\mathcal A[\tilde f_x](s)\big), \]
which gives us with the representation~\eqref{eq:inverse_abel} of the inverse Abel transform that
\begin{equation}\label{eq:spmeans_and_2d_wave_equation}
\mathcal M_2[f](x;r) = \frac1{\pi r}\mathcal A[\tilde p_x](\tfrac1r)=\frac2\pi\int_0^r\frac{p(x,t)}{\sqrt{r^2-t^2}}\d t,\quad x\in\R^2,\;r>0,
\end{equation}
where $\tilde p_x(\tfrac1t) = t^2p(x;t)$; and for $n=3$, we find
\[ \mathcal M_3[f](x;r) = \frac1r\int_0^rp(x;t)\d t,\quad x\in\R^3,\;r>0. \]

\subsection{The Radon Transform}
The Radon transform for a function $g:\R^2\to\R$ is defined by
\begin{equation}\label{eq:radon_transform}
\mathcal R[g]:\R\times S^1\to\R,\quad\mathcal R[g](r,\theta) = \int_{-\infty}^\infty g(r\theta+u\theta^\perp)\d u.
\end{equation}

The Radon transform can be inverted and we have the explicit formula
\[ \mathcal R^{-1}[G](\xi) = -\frac1{(2\pi)^2}\int_{-\infty}^\infty\frac1r\int_{S^1}\partial_rG(r+\left<\theta,\xi\right>,\theta)\d\theta\d r \]
for the inverse Radon transform, see e.g.~\cite{Nat01}.

\subsection{Mathieu Functions}\label{sec:Mathieu_equation}
The Mathieu Functions are solutions of the Mathieu equation
\begin{equation}\label{eq:Mathieu}
u''(s)+(a-2q\cos(2s))u(s)=0.
\end{equation}
However, we are only interested in $2\pi$-periodic solutions. It is known, see e.g.~\cite{AbrSte64}, that for a fixed value $q\ge0$, there only exists a $2\pi$-periodic solution of the equation~\eqref{eq:Mathieu} for a discrete set of values $a\in\R$. Conventionally, the values $a$ for which an even $2\pi$-periodic solution exists, are labeled in increasing order as $a_n(q)$ with the Mathieu cosine functions $s\mapsto\ce_n(s;q)$ as corresponding solutions, $n\in\N_0$; and the values $a$ for which we have an odd $2\pi$-periodic solution are (again in increasing order) called $b_n(q)$ with the Mathieu sine functions $s\mapsto\se_n(s,q)$ as corresponding solutions, $n\in\N$. The normalization of the solutions $\ce_n$ and $\se_n$ is chosen to be
\[ \int_0^{2\pi}\ce_n(s;q)^2\d s = \pi\quad\text{and}\quad\int_0^{2\pi}\se_n(s;q)^2\d s = \pi. \]
Thus, since the functions are the eigenfunctions of the symmetric operator $\partial_{ss}+2q\cos(2s)$, the functions $\frac1{\sqrt\pi}\ce_n(\cdot;q)$, $\frac1{\sqrt\pi}\se_{n+1}(\cdot;q)$, $n\in\N_0$, form for every $q\ge0$ a complete orthonormal system of $L^2([0,2\pi])$.

\bibliography{citations}
\bibliographystyle{abbrv}
\end{document}

%% file: point_of_tangency.tex
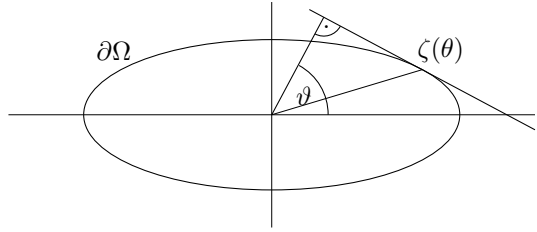
\begin{figure}[ht]
\begin{center}
\begin{tikzpicture}[scale=0.5]
% ellipse
\draw(0,0) circle (5cm and 2cm);
\node at (-4.2,1.7) {$\partial\Omega$};
%\node at (-1.5,0.7) {$\Omega$};
% x axis
\draw[very thin](-7,0) -- (7,0);
\draw[very thin](0,-3) -- (0,3);
% line to point of tangency
\draw[very thin](0,0) -- (4,1.2);
\node at (4.5,1.7) {$\zeta(\theta)$};
% tangent
\draw(1,2.8) -- (7,-0.4);
% normal to tangent
\draw(0,0) -- (1.384,2.594);
% label orthogonal angle
%\draw(1.384,2.594)+(-118.072:0.5) arc (-118.072:-28.072:0.5);
\draw(1.384,2.594)+(-118:0.5) arc (-118:-28:0.5);
\node at (1.457,2.355) {.};
%(1.384,2.594)+(-73.072:0.25) {$\cdot$};
% label angle \theta
\node at (31:1) {$\vartheta$};
\draw(1.5,0) arc (0:61:1.5);
%\draw(1.5,0) arc (0:61.918:1.5);
\end{tikzpicture}
\end{center}
\caption{Definition of the point $\zeta(\theta)$, $\theta=(\cos\vartheta,\sin\vartheta)$.}\label{fig:zeta}
\end{figure}